\documentclass[a4paper,12pt]{amsart}

\title[Non-doubling measures] {\bf Extended multifractal formalism of some
  non-doubling measures} \author{Shuang Shen} \date{December 2014}

\usepackage{latexsym,amsmath,amsfonts,mathrsfs}
\usepackage{amsthm}
\usepackage{amssymb}
\usepackage{amstext}
\usepackage{cite}
\usepackage{amsopn}
\DeclareMathOperator{\Dim}{\mathrm{Dim}}
\DeclareMathOperator*{\esssup}{\mathrm{ess\,sup}}
\newcommand{\ball}{\mathrm{B}}
\newtheorem{thm}{Theorem}
\newtheorem{lemma}[thm]{Lemma}
\newtheorem{cor}[thm]{Corollary}
\theoremstyle{remark}
\newtheorem{remark}[thm]{Remark}
\begin{document}

\thanks{The research is supported by the National Natural Science Foundation 
of China, No. 11271223, No. 11431007 and No. 11201256}

\address{Department of Mathematical Sciences, Tsinghua University, Beijing 100084, China}
\email{shens10@mails.tsinghua.edu.cn}

\begin{abstract}  

In a previous work~\cite{She} we constructed measures on symbolic
spaces which satisfy an extended multifractal formalism (in the sense 
that Olsen's functions~$b$ and~$B$ differ and that their Legendre
transforms have the expected interpretation in terms of
dimensions). These measures are composed with a Gray code and
projected onto the unit interval so to get doubling measures. Then we
were able to show that the projected measure has the same Olsen's
functions as the one it comes from and that it also fulfills the
extended multifractal formalism. Here we show that the use of a Gray
code is not necessary to get these results, although dealing with non
doubling measures.

\end{abstract}

\maketitle

\section{Introduction}

Ben Nasr, Bhouri, and Heurteaux~\cite{Ben} constructed a class of
measures whose Olsen's~$b$ and~$B$ functions differ. They first
consider a Markov measure on the symbolic space $\{0,1\}^{\mathbb N}$
and project it on $[0,1]$ by using the usual dyadic
representation~$\gamma$ of numbers. The Markov rules are chosen so
that the projected measure is doubling. Actually this is equivalent to
the following construction.

Given two different numbers $a, a'\in (0,1)$ and an increasing
sequence $T_k$ of positive integers such that $\lim_{k\rightarrow
  \infty} T_{k+1}/T_k = \infty$, consider the measure $\mu$ on the
symbolic space so defined: for any $w=w_1\cdots w_n\in \{0,1\}^n$,
\begin{equation*}
 \mu([w]) = \prod_{j=1}^n \Bigl(p_j^{1-w_j}(1-p_j)^{w_j}\Bigr),
\end{equation*}
where $[w]$ stands for the cylinder defined by $w$, and, $p_j=a$ if
$T_{2k-1}\leq j< T_{2k}$, and $p_j=a'$ if $T_{2k}\leq j< T_{2k+1}$ for
some~$k$. If we denote by~$g$ the Gray code (this is a map from the
symbolic space into itself which allows to enumerate cylinders of the
same size in such a way that one passes from an element to the next
one by flipping one digit only), then the measure~$\nu$ considered
in~\cite{Ben} is the image of~$\mu$ under~$\gamma\circ g$.

Since~$\nu$ is doubling, when concerning multifractal analysis, it inherits
from~$\mu$ (see~\cite{Ben,Ben2,She}).

On the other hand, Barral, Ben Nasr and Peyri{\`e}re~\cite{Bar} 
showed that the projection under~$\gamma$ of a Bernoulli 
measure satisfies the multifractal formalism although it is not 
doubling. This is why we investigate whether it was necessary 
to use a Gray code to obtain a measure on~$[0,1]$ satisfying an 
extended multifractal formalism (or refined multifractal formalism, 
according to Barral's terminology \cite{Bar2}). Indeed, we are 
going to show that projections under~$\gamma$ of inhomogeneous
multinomial measures have this property.

Barral \cite{Bar2} proved that, given two convex functions fulfilling 
fairly general conditions, there exists a compactly supported 
(always supported on a Cantor set), positive, and finite Borel measure 
$\rho$ on $\mathbb{R}$ whose $\tau_\rho$ and $\underline{\tau}_\rho$ 
functions are just the two given functions. Also, the author mentioned 
that for the measure $\rho$ possessing the weak doubling properties 
(see Inequality~\eqref{key}), if $\dim X_\rho(\alpha)=\underline{\tau}_
\rho^\ast(\alpha)$ for all $\alpha$ over its domain, then 
$b_\rho=\underline{\tau}_\rho$; similarly if $\Dim X_\rho(\alpha)=
\tau_\rho^\ast(\alpha)$ for all $\alpha$ over its domain, then 
$B_\rho=\tau_\rho$ (these notations will be reminded later whereas for the 
definition of domain, the reader is referred to \cite{Bar2}).

In contrast, for the measures we consider, generally the $\underline{\tau}$ 
functions are not convex (see Theorem \ref{t1}), and the condition 
$\Dim X(\alpha)= \tau^\ast(\alpha)$ does not always hold (see 
Theorem \ref{t2}). Moreover, the measures we construct in this article 
have full support.

This article is organized as follows. In Section~\ref{sec2}, we 
recall the basic notations, definitions, and the constructions of 
inhomogeneous multinomial measures. In Section~\ref{sec3}, 
we present our main results. Then we prove the two results in 
Section~\ref{sec4} and Section~\ref{sec5} respectively.

\section{Recollections: Notations and definitions}\label{sec2}
\subsection{The Olsen's measures and functions}

We work on a metric space $(\mathbb{X},d)$ possessing the Besicovitch
property:
\smallskip

\textsl{There exists a constant $C_\ball\in \mathbb{N}$ such that,
  given any bounded subset $\{x_i\}_{i\in I}\subseteq \mathbb{X}$ and
  any collection $\{\ball(x_i,r_i)\}_{i\in I}$ of balls in
  $\mathbb{X}$, one can extract from it $C_\ball$ countable families
  $\{\{\ball(x_{j,k},r_{j,k})\}_{k\geq 1}\}_{1\leq j\leq
    C_\ball}$ so that
\begin{enumerate}
\item[--] $ \bigcup_{j,k}\ball(x_{j,k},r_{j,k})\supseteq \{x_i\}_{i\in I}$,
\item[--] for any $j$ and $k\neq k'$, $\ball(x_{j,k},r_{j,k})\cap
  \ball(x_{j,k'},r_{j,k'})=\emptyset$.
\end{enumerate}}
\bigskip

It is known that Euclidean spaces and ultrametric spaces fulfill this
condition.

Let $\mu$ be a Borel probability measure on $\mathbb{X}$. Denote 
by $S_\mu$ the support of the measure $\mu$. For any 
$\alpha\in\mathbb{R}$, we denote by $X_\mu(\alpha)$ the level set 
of points whose local H\"older exponents assume the value~$\alpha$:
$$X_\mu(\alpha)=\left\{x\in S_\mu:
\limsup_{r\rightarrow 0}\frac{\log\mu(\ball(x,r))}{\log r}=
\liminf_{r\rightarrow 0}\frac{\log\mu(\ball(x,r))}{\log r}=\alpha\right\}$$

For $q, t\in \mathbb{R}$ and $\delta >0$, we shall use
the measures and premeasures $\overline{\mathscr{H}}_{\mu,\delta}^{q,t}$,
$\overline{\mathscr{H}}_{\mu}^{q,t}$, $\mathscr{H}_{\mu}^{q,t}$
and $\overline{\mathscr{P}}_{\mu,\delta}^{q,t}$,
$\overline{\mathscr{P}}_{\mu}^{q,t}$, $\mathscr{P}_{\mu}^{q,t}$
introduced by Olsen~\cite{Ols} and whose definitions are recalled in~\cite{She}.

The functions $\mathscr{H}_{\mu}^{q,t}$, $\mathscr{P}_{\mu}^{q,t}$ and
$\overline{\mathscr{P}}_{\mu}^{q,t}$ provide each subset $E$ of
$\mathbb{X}$ with dimensional indices:
$$
\begin{array}{rccclll}
\medskip
b_{\mu,E}(q) &=& \sup\{s\ :\ \mathscr{H}_\mu^{q,s}(E)=
\infty\} &=& \inf\{s\ :\ \mathscr{H}_\mu^{q,s}(E)=0\},\\
\medskip
B_{\mu,E}(q) &=& \sup\{s\ :\ \mathscr{P}_\mu^{q,s}(E)=
\infty\} &=& \inf\{s\ :\ \mathscr{P}_\mu^{q,s}(E)=0\},\\
\tau_{\mu,E}(q) &=& \sup\{s\ :\ \overline{\mathscr{P}}_\mu^{q,s}(E)=
\infty\} &=& \inf\{s\ :\ \overline{\mathscr{P}}_\mu^{q,s}(E)=0\}.
\end{array}
$$

One sees that the Olsen's function $b_{\mu,E}(q)$ is a multifractal extension 
of the Hausdorff dimension $\dim E$ whereas $B_{\mu,E}(q)$ is a multifractal 
extension of the packing dimension $\Dim E$. For simplicity, when the context 
is clear, we will write $b_\mu=b_{\mu,S_\mu}$,
$B_\mu=B_{\mu,S_\mu}$, and $\tau_\mu=\tau_{\mu,S_\mu}$. Some basic
properties can be found in \cite{Ols}.

There is an alternate way to compute the function $\tau_\mu$ according to
\cite{Ben2,Pey}. Fix $\lambda<1$, one has
\begin{multline}\label{alt}
\tau_{\mu}(q)=\limsup_{\delta\rightarrow
0}\frac{-1}{\log \delta}\log\sup\left\{\sum_i
\mu(\ball_i)^q\ :\  \right .\\
 \left .(\ball_i)_i \textrm{ is a packing of }S_\mu \textrm{ with }
 \lambda\delta<r_i\leq\delta\vphantom {\sum_i}\right\},
\end{multline}
where $r_i$ is the radius of the ball $\ball_i$.

\subsection{Inhomogeneous multinomial measures}
\subsubsection{The mixed symbolic spaces}
Let $c_1,c_2\geq 2$ be two positive integers,
$\mathscr{A}_1=\{0,1,\cdots,c_1-1\}$, $\mathscr{A}_2=\{0,1,\cdots,c_2-1\}$
be two alphabets. Fix in this article a sequence of integers $(T_k)$ such that
$$T_1=1, \,T_k<T_{k+1} \textrm{ and } \lim_{k\rightarrow\infty}
T_{k+1}/T_k=+\infty.$$

Consider the set of infinite words
$$\partial\mathscr{A}_{1,2}^\ast = 
\mathscr{A}_1^{T_2-T_1}\times \mathscr{A}_2^{T_3-T_2}\times
\mathscr{A}_1^{T_4-T_3}\times\cdots =\prod_j Y_j,$$ where
\begin{enumerate}
\item[--] if $T_{2k-1}\leq j< T_{2k}$ for some $k$, \,$Y_j=\mathscr{A}_1$,
\item[--] if $T_{2k}\leq j< T_{2k+1}$ for some $k$, \,$Y_j=\mathscr{A}_2$.
\end{enumerate}

We call $\partial\mathscr{A}_{1,2}^\ast$ the mixed symbolic
space with respect to the triplet $\{\mathscr{A}_1,\mathscr{A}_2,(T_k)\}$. 
Also we denote by $\mathscr{A}_{1,2}^n$ the set of words of length~$n$
(by convention the empty word $\epsilon$ has length~$0$), 
and $\mathscr{A}_{1,2}^\ast$ the set of finite words, i.e. 
\begin{eqnarray*}
\mathscr{A}_{1,2}^0 &=& \{\epsilon\},\\
\mathscr{A}_{1,2}^n &=& \prod_{j=1}^n Y_j,  \; n\geq 1,\\
\mathscr{A}_{1,2}^\ast &=& \bigcup_{n\geq 0}\mathscr{A}_{1,2}^n.
\end{eqnarray*}  

The length of a finite word $w$ is denoted 
by $\ell(w)$. If $w=\varepsilon_1\cdots\varepsilon_k\cdots$ is an infinite 
word or a finite word of length larger than $k$, the $k$-prefix of~$w$ 
is denoted by $w|_k= \varepsilon_1\cdots\varepsilon_k$. And for any 
word $w$, by convention one has $w|_0=\epsilon$, where $\epsilon$ is 
the empty word. If $w$ and $v$ are two words, $w\wedge v$ stands for 
their largest common prefix. 

Let $N_n$ be the number of integers $j\leq n$ such that
$Y_j=\mathscr{A}_1$. We can immediately get that
\begin{equation}\label{time}
\liminf_{n\rightarrow\infty}\frac{N_n}{n}=0 \;\textrm{  and  } \;
\limsup_{n\rightarrow\infty}\frac{N_n}{n}=1.
\end{equation}

For any two different elements $w,v\in \partial\mathscr{A}_{1,2}^\ast$
with $\ell(w\wedge v)=n$, we define $d(w,v)=c_1^{-N_n}c_2^{-(n-N_n)}$. 
It is easy to check that this defines an
ultrametric distance on~$\partial\mathscr{A}_{1,2}^\ast$.

Each finite word $w\in\mathscr{A}_{1,2}^\ast$ defines a cylinder
$[w]=\{x\in\partial\mathscr{A}_{1,2}^\ast:x|_{\ell(w)}=w\}$, which can also be
viewed as a ball; and the diameter of this ball is denoted by $|w|$.
For a Borel measure $\mu$ on $\partial\mathscr{A}_{1,2}^\ast$, we simply
write $\mu([w])=\mu(w)$. Thus we identify the Borel measure $\mu$ on
$\partial\mathscr{A}_{1,2}^\ast$ with a mapping from $\mathscr{A}_{1,2}^\ast$ 
to $[0,+\infty]$ subject to the following compatibility condition 
$$\mu(w)=\sum_{\substack{x\in\mathscr{A}_{1,2}^{n+1} \\ x|_n=w}}\mu(x),\; 
\textrm{ for any } n\geq 0 \textrm{ and } w\in\mathscr{A}_{1,2}^n.$$

One sees that when the alphabets $\mathscr{A}_1=\mathscr{A}_2=\mathscr{A}$, 
the mixed symbolic space $\partial\mathscr{A}_{1,2}^\ast$ becomes 
ordinary symbolic space $\partial\mathscr{A}^\ast$.

Since the radii of balls are discrete on the mixed symbolic space, the
computation of~$\tau_\mu$ according to Formula~\eqref{alt} is quite
easy.  Indeed, for any $n$, we take any element
$w\in\mathscr{A}_{1,2}^n$ and define $\tau_{\mu,n}$ by the following
formula
\begin{equation}\label{tau}
\sum_{z\in\mathscr{A}_{1,2}^n} \mu(z)^q=|w|^{-\tau_{\mu,n}(q)}.
\end{equation}
Then
$$\tau_\mu(q)=\limsup_{n\rightarrow\infty} \tau_{\mu,n}(q).$$
Also, we denote
$$\underline{\tau}_\mu(q)=\liminf_{n\rightarrow\infty}
\tau_{\mu,n}(q).$$

\subsubsection{Image measures}
There is a natural map from $\partial\mathscr{A}_{1,2}^\ast$ onto
$\mathbb{R}$. Consider the map $\gamma$ which sends the element 
$x=\varepsilon_1 \varepsilon_2\cdots \varepsilon_n\cdots $ to the real
number $\sum_{n\geq 1} \varepsilon_n c_1^{-N_n}c_2^{-(n-N_n)}$. 
This map sends $n$-cylinders to basic intervals of $n$-th generation and 
thus defines a function on $\mathscr{A}_{1,2}^\ast$, still denoted by $\gamma$. 
To be precise, one can assign to each $w\in \mathscr{A}_{1,2}^n$
an integer $\iota(w)$ such that
$$\gamma(w)=\left[\iota(w)c_1^{-N_n}c_2^{-(n-N_n)},
(\iota(w)+1)c_1^{-N_n}c_2^{-(n-N_n)}\right].$$

Now if $\mu$ is a Borel probability measure on
$(\partial\mathscr{A}_{1,2}^\ast,d)$, we have its projection $\nu$ on
$(\mathbb{R},|\cdot|)$. This is the image measure
$\nu=\gamma_\ast(\mu)$, defined by $\nu(E)=\mu(\gamma^{-1}(E))$ for
any Borel set $E\subseteq [0,1]$. 

In particular case, the projection under~$\gamma$ of a Bernoulli
measure satisfies the multifractal formalism, although it is not
doubling. In fact, J.~Barral, F.~Ben Nasr and J.~Peyri{\`e}re proved
the following stronger result.

\begin{thm}[see \cite{Bar}]\label{t0}
Let $\mu$ be a continuous quasi-Bernoulli measure on
$\partial\mathscr{A}^\ast$. Then both measures $\mu$ and
$\nu=\gamma_\ast(\mu)$ obey the multifractal formalism everywhere and
one has
$$b_\nu=B_\nu=\tau_\nu=b_\mu=B_\mu=\tau_\mu.$$
\end{thm}
\bigskip

However in general, the measure we are going to study is not
quasi-Bernoulli.

\subsubsection{Inhomogeneous multinomial measures}
Given two groups of real numbers $a_i, b_j\in (0,1)$($i=1,
\cdots,c_1,\,j=1,\cdots,c_2$) satisfying
$$a_1+\cdots+a_{c_1}=b_1+\cdots+b_{c_2}=1,$$
we define a probability measure~$\mu$ on $\partial\mathscr{A}_{1,2}^\ast$ 
that we call an inhomogeneous multinomial measure as explained below.  
As in~\cite{She}, for every cylinder $[\varepsilon_1
\varepsilon_2\cdots\varepsilon_n]$, we set
$$\mu(\varepsilon_1\cdots\varepsilon_n)=\prod_{j=1}^n p_{j},$$
where
\begin{enumerate}
\item[--] if $T_{2k-1}\leq j< T_{2k}$ for some $k$, \,$p_j=a_{\varepsilon_j+1}$,
\item[--] if $T_{2k}\leq j< T_{2k+1}$ for some $k$, \,$p_j=b_{\varepsilon_j+1}$.
\end{enumerate}
\medskip

Then we compute the $\tau_\mu$ function.  As previously, let $N_n$
stand for the number of integers $j\leq n$ such that $p_j\in
\{a_1,\cdots,a_{c_1}\}$, then we obtain by Formula~\eqref{tau} that
$$\tau_{\mu,n}(q)=\frac{\frac{N_n}{n}\log(a_1^q+\cdots+a_{c_1}^q)}
{\frac{N_n}{n}\log c_1+(1-\frac{N_n}{n})\log c_2}+
\frac{(1-\frac{N_n}{n})\log(b_1^q+\cdots+b_{c_2}^q)}
{\frac{N_n}{n}\log c_1+(1-\frac{N_n}{n})\log c_2}.$$
 This combined with \eqref{time} implies
$$\tau_\mu(q)=\max\{\log_{c_1}(a_1^q+\cdots+a_{c_1}^q),\log_{c_2}(b_1^q+
\cdots+b_{c_2}^q)\},$$
$$\underline{\tau}_\mu(q)=\min\{\log_{c_1}(a_1^q+\cdots+a_{c_1}^q),
\log_{c_2}(b_1^q+\cdots+b_{c_2}^q)\}.$$

Recall the following fact.

\begin{thm}[see \cite{She}]\label{t1} 
One has
$$B_\mu(q)=\tau_\mu(q)=\max\{\log_{c_1}(a_1^q+\cdots+a_{c_1}^q),\,
\log_{c_2}(b_1^q+\cdots+b_{c_2}^q)\},$$
$$b_\mu(q)=\underline{\tau}_\mu(q)=\min\{\log_{c_1}(a_1^q+\cdots+a_{c_1}^q),\,
\log_{c_2}(b_1^q+\cdots+b_{c_2}^q)\}.$$
\end{thm}
\bigskip

We need the following auxiliary measures.
\begin{lemma}[see \cite{Ben,She}]\label{key l1}
For any $q\in\mathbb{R}$, there is a probability
measure $\mu_q$ on $\partial\mathscr{A}_{1,2}^\ast$ and a
subsequence of integers $(n_k)_{k\geq 1}$, such that
$$\mu_q(w)= \mu(w)^q |w|^{\tau_{\mu,n}(q)},\textrm{ if } w\in
  \mathscr{A}_{1,2}^n.$$
Moreover,
$$\mu_q(w)\leq \mu(w)^q |w|^{\underline{\tau}_\mu(q)},\textrm{ if } w\in
  \mathscr{A}_{1,2}^n,$$
and for every $\varepsilon> 0$,
$$\mu_q(w)\leq \mu(w)^q |w|^{\tau_\mu(q)-\varepsilon}, \textrm{ if } w\in
\mathscr{A}_{1,2}^{n_k} \textrm{ with $k$ large}.$$
\end{lemma}
\medskip

We also point out that all these measures $\mu, \mu_q$ are continuous.

\section{Main results}\label{sec3}

Let us state our main results. Denote 
$$\theta_a(q)=\log_{c_1}(a_1^q+\cdots+a_{c_1}^q),$$
$$\theta_b(q)=\log_{c_2}(b_1^q+\cdots+b_{c_2}^q).$$
And
$$s_1=\min \left\{\max_{1\leq i\leq c_1}\{a_i\},\,\max_{1\leq j\leq
  c_2}\{b_j\}\right\},$$
$$s_2=\max \left\{\min_{1\leq i\leq c_1}\{a_i\},\,\min_{1\leq j\leq
  c_2}\{b_j\}\right\}.$$
Let $f^\ast(x)= \inf_y (xy+f(y))$ denote the Legendre transform of the
function $f$.

\begin{thm}\label{t2} 
Let $\mathscr{A}_1=\{0,1,\cdots,c_1-1\}$, $\mathscr{A}_2=\{0,1,\cdots,c_2-1\}$ 
and let $\mu$ be the probability measure on $\partial\mathscr{A}_{1,2}^\ast$ 
taken from Theorem~\ref{t1}. Denote $\nu=\gamma_\ast(\mu)$. Then for every
$q\in \mathbb{R}$,
$$B_\nu(q)=B_\mu(q)=\max\{\theta_a(q),\,\theta_b(q)\},$$
$$b_\nu(q)=b_\mu(q)=\min\{\theta_a(q),\,\theta_b(q)\}.$$
\end{thm}

\begin{thm}\label{t3}
For any $\alpha\in(-\log s_1,-\log s_2)$, we have
$$\dim X_\mu(\alpha)=\dim X_\nu(\alpha)=b_\mu^\ast(\alpha)=b_\nu^\ast(\alpha).$$
And for $\alpha\in(-\log s_1,-\log s_2)$ subject to
$$\max\{\theta^\ast_a(\alpha),\,\theta^\ast_b(\alpha)\}=B^\ast_\mu(\alpha),$$
we have
$$\Dim X_\mu(\alpha)=\Dim X_\nu(\alpha)=B_\mu^\ast(\alpha)=B_\nu^\ast(\alpha).$$
\end{thm}

\section{Proof of Theorem \ref{t2}}\label{sec4}

To avoid tedious notations, we write the proof with $c_1=c_2=2$. The
reader will realize that the general case can be handled with minor
modifications. For $n\geq 0$, we denote by $\mathscr{F}_n$ the family
of basic intervals of $n$-th generation:
$$\mathscr{F}_n=\left\{\left[\frac{j}{2^n},\frac{j+1}{2^n}\right],\,0\leq
j\leq 2^n-1\right\}.$$

For $x\in (0,1)$, denote by $I_n(x)$ the basic interval of $n$-th
  generation which contains $x$ when $x$ is not of the form
  $k2^{-m}$. When $x$ is dyadic there are two basic intervals of
  $n$-th generation containing $x$ and we choose the left one to be
  $I_n(x)$. We also denote by $I_n(x)^\pm$ the basic intervals of the
  same generation, which are the right and left neighbors of $I_n(x)$
  respectively. 

The following key lemma says that the measure $\nu=\gamma_\ast(\mu)$
exhibits some weak doubling behaviors.

\begin{lemma}\label{key l2} 
Let $\varpi$ be a continuous probability measure on $[0,1]$ satisfying that, 
there exists a positive constant $C_1<1$, such that for any $n\geq 0$, 
for any $J\in\mathscr{F}_n$, and for $I\in\mathscr{F}_{n+1}$ with 
$I\subseteq J$, one has $\varpi(I)\leq C_1 \varpi(J)$.

Then there exists a constant $C_0>1$ 
such that for $\varpi$-almost every $x\in [0,1]$, when $n$ is large enough,
\begin{equation}\label{key}
C_0^{-(\sqrt{n}+1)}\leq \frac{\nu(I_n(x))}{\nu(I_n(x)^\pm)}\leq
C_0^{\sqrt{n}+1}.
\end{equation}
\end{lemma}

\proof
Denote $\langle n\rangle=\lfloor n-\sqrt{n}\rfloor.$
We define the set
\begin{multline*}
E_n=\{x\in [0,1]:\,\ell(\gamma^{-1}(I_n(x))\wedge\gamma^{-1}(I_n(x)^+))< 
\langle n\rangle \textrm{ or } \\
\ell(\gamma^{-1}(I_n(x))\wedge\gamma^{-1}(I_n(x)^-))< \langle n\rangle\}.
\end{multline*}

We fix $n$ for the moment and denote the basic intervals of 
$\langle n\rangle$-th generation, from left to right, by 
$I_1, \cdots, I_{k_n}$, where $k_n=2^{\langle n\rangle}$. The 
leftmost and rightmost basic subinterval, of $n$-th generation, of $I_j$ are 
denoted by $I_j^{(1)}$ and $I_j^{(2)}$ respectively. So by definition
$$E_n\subseteq\bigcup_{j=1}^{k_n} I_j^{(1)}\cup I_j^{(2)}.$$

As an example, we compare the $\varpi$-mass of $I_j^{(1)}$ with its
mother interval, denoted by $J$. Note
that $I_j^{(1)}\in \mathscr{F}_n$ and $J\in \mathscr{F}_{n-1}$. Thus
it follows that $\varpi(I_j^{(1)})\leq C_1 \varpi(J)$,
which implies $\varpi(I_j^{(1)})\leq C_1^{\sqrt{n}}\varpi(I_j)$. Moreover,
$$\varpi(E_n)\leq\sum_{j=1}^{k_n}
(\varpi(I_j^{(1)})+\varpi(I_j^{(2)}))\leq 2C_1^{\sqrt{n}}
\sum_{j=1}^{k_n} \varpi(I_j)=2C_1^{\sqrt{n}},$$ 
from which it follows
$$\sum_{n\geq 1}\varpi(E_n)\leq \sum_{n\geq 1}2C_1^{\sqrt{n}}<+\infty.$$
By Borel-Cantelli lemma, one immediately gets
$$\varpi(\limsup_n E_n)=0,$$
which implies
$$\varpi(\liminf_n E_n^{\mathsf c})=1,$$
where $E_n^{\mathsf c}=[0,1]\setminus E_n$.

To conclude, one recalls that
$\ell(\gamma^{-1}(I_n(x))\wedge\gamma^{-1}(I_n(x)^\pm))\geq \langle
n\rangle$ if $x\in E_n^{\mathsf c}$, meaning that the $\nu$-masses will not
differ too much. In fact, denote
$$C_0=\max\left\{\frac{a_1}{a_2},\frac{a_2}{a_1},\frac{b_1}{b_2},
\frac{b_2}{b_1}\right\}.$$
Then $C_0>1$ and
$$C_0^{-(\sqrt{n}+1)}\leq \frac{\nu(I_n(x))}{\nu(I_n(x)^\pm)}\leq C_0^{\sqrt{n}+1}.$$
\qed

\begin{cor}\label{coro}
Recall the auxiliary measures $\mu_q$ presented in Lemma \ref{key l1}. 
For any $q\in\mathbb{R}$, let $\nu_q = \gamma_\ast(\mu_q)$, then for 
$\nu_q$-almost every $x\in [0,1]$, when $n$ is large enough, 
Inequality~\eqref{key} holds.
\end{cor}

\proof
For any $q$, let
$$C_1(q)=\max\left\{\frac{a_1^q}{a_1^q+a_2^q},\frac{a_2^q}
{a_1^q+a_2^q},\frac{b_1^q}{b_1^q+b_2^q},\frac{b_2^q}{b_1^q+b_2^q}\right\}.$$

Then $C_1(q)<1$ and for any $J\in\mathscr{F}_n$, for $I\in\mathscr{F}_{n+1}$ with 
$I\subseteq J$, one obtains by Lemma \ref{key l1} that $\nu_q(I)\leq C_1(q) \nu_q(J)$.
\qed

\subsection{$b_\nu(q)=\underline{\tau}_\mu(q)$}
We first consider the easier part $b_\nu(q)\leq \underline{\tau}_\mu(q)$. 
For any $\varepsilon>0$, choose a subsequence $\{n_k\}$ such that 
$\tau_{\mu,n_k}(q)<\underline{\tau}_\mu(q)+\varepsilon$, for every $k\geq
1$. Take any subset $F\subseteq S_\nu=[0,1]$, and we choose the natural
centered $2^{-n_k}$-covering of $F$, which is a set of all basic intervals 
of ${n_k}$-th generation. Now
$$\overline{\mathscr{H}}^{q,\underline{\tau}_\mu(q)+\varepsilon}_{\nu,2^{-n_k}}(F)
\leq\sum_{I \in\mathscr{F}_{n_k}}\nu(I)^q 2^{-n_k(\underline{\tau}_\mu(q)+\varepsilon)}=2^{n_k(\tau_{\mu,n_k}(q)-\underline{\tau}_\mu(q)-\varepsilon)}\leq 1,$$
which implies
$$\overline{\mathscr{H}}^{q,\underline{\tau}_\mu(q)+\varepsilon}_\nu(F)\leq 1,$$
and
$$\mathscr{H}^{q,\underline{\tau}_\mu(q)+\varepsilon}_\nu(S_\nu)\leq 1.$$
Thus $b_\nu(q)\leq \underline{\tau}_\mu(q)$ since $\varepsilon$ is arbitrary.

Next we turn to the opposite part $b_\nu(q)\geq \underline{\tau}_\mu(q)$.

\begin{lemma}\label{l3}
Let $q\in \mathbb{R}$. For any $\varepsilon>0$, and for $\nu_q$-almost every $x$, 
when $r$ is small enough,
$$\nu_q(\ball(x,r))\leq \nu(\ball(x,r))^q r^{\underline{\tau}_\mu(q)-\varepsilon}.$$
Then we have $b_\nu(q)\geq\underline{\tau}_\mu(q)$.
\end{lemma}

\proof
For $\nu_q$-almost every $x\in [0,1]$, when $r$ is small enough, we can 
choose $n$ large enough such that $2^{-(n+2)}<r\leq 2^{-(n+1)}$. So there 
exist at most two basic intervals $I_1, I_2\in\mathscr {F}_n$ such that 
$\ball(x,r)\subseteq I_1\cup I_2$; on the other hand, $\ball(x,r)$ must contain 
at least one basic interval of $(n+2)$-th generation, denoted by $I_3$. 
It is no restriction to assume that $x\in I_1$ and $x\in I_3$. Note that by 
Corollary \ref{coro}, $\nu(I_1)$ and $\nu(I_2)$ do not differ too much since 
$n$ is large. With the help of Lemma \ref{key l1}, we have
\begin{eqnarray*}
\nu_q(\ball(x,r))\leq \nu_q(I_1)+\nu_q(I_2)&\leq&\nu(I_1)^q
|I_1|^{\underline{\tau}_\mu(q)}+\nu(I_2)^q
|I_2|^{\underline{\tau}_\mu(q)}\\
&=&(\nu(I_1)^q+\nu(I_2)^q)|I_1|^{\underline{\tau}_\mu(q)},
\end{eqnarray*}
where $|I|$ stands for the length of the interval $I$.

When $q<0$, we first consider the case where $\nu(I_1)\geq
\frac{1}{2}\nu(\ball(x,r))$. We then have
$$\nu(I_2)\geq \frac{1}{C_0^{\sqrt{n}+1}} \nu(I_1)\geq 
\frac{1}{2C_0^{\sqrt{n}+1}}\nu(\ball(x,r)),$$
which implies
\begin{eqnarray*}
\nu_q(\ball(x,r))&\leq& \left(\frac{1}{2^q}+\frac{1}{2^q C_0^{q(\sqrt{n}+1)}}\right)
\nu(\ball(x,r))^q |I_1|^{\underline{\tau}_\mu(q)}\\
&\leq& C(q)C_0^{-q\sqrt{n}}\nu(\ball(x,r))^q r^{\underline{\tau}_\mu(q)}\\
&=&\nu(\ball(x,r))^q r^{\underline{\tau}_\mu(q)-\varsigma(q,r)},
\end{eqnarray*}
where
$$C(q)=2^{-(q-1)}C_0^{-q}\max\{2^{\underline{\tau}_\mu(q)},
4^{\underline{\tau}_\mu(q)}\},$$
and
$$\varsigma(q,r)=\frac{\log C(q)-q\sqrt{n}\log C_0}{-\log r}.$$

So for any $\varepsilon>0$, take $r$ small enough such that 
$\varsigma(q,r)<\varepsilon$, then we have
$$\nu_q(\ball(x,r))\leq \nu(\ball(x,r))^q r^{\underline{\tau}_\mu(q)-\varepsilon}.$$

With the same method, in the case where $\nu(I_2)\geq
\frac{1}{2}\nu(\ball(x,r))$, we can get the same results. 

When $q\geq 0$, we use the fact that $\nu(\ball(x,r))\geq \nu(I_3)$. 
Since $I_1$ is the grandmother interval of $I_3$, we have 
$\nu(I_3)\geq C_2^2\nu(I_1)$ (where $C_2=\min\{a_1,a_2,b_1,b_2\}<1$), 
which implies
$$\nu(I_1)\leq C_2^{-2}\nu(I_3)\leq C_2^{-2}\nu(\ball(x,r)),$$
and thus
$$\nu(I_2)\leq C_0^{\sqrt{n}+1}\nu(I_1)\leq C_0^{\sqrt{n}+1}C_2^{-2}\nu(\ball(x,r)).$$

Then in the very same way as above, we also conclude that for any 
$\varepsilon>0$, when $r$ is small enough, we have
$$\nu_q(\ball(x,r))\leq \nu(\ball(x,r))^q r^{\underline{\tau}_\mu(q)-\varepsilon}.$$
\medskip

Finally, we have to show $b_\nu(q)\geq\underline{\tau}_\mu(q)$. 
For any $\varepsilon>0$, denote
$$W=\{x\in [0,1]:\,\exists r_x>0, \forall r<r_x, \nu_q(\ball(x,r))\leq 
\nu(\ball(x,r))^q r^{\underline{\tau}_\mu(q)-\varepsilon}\},$$
$$W_n=\{x\in W:\,\forall r<1/n, \nu_q(\ball(x,r))\leq 
\nu(\ball(x,r))^q r^{\underline{\tau}_\mu(q)-\varepsilon}\},\,\forall n\in \mathbb{N}.$$
Then $\nu_q(W)=1$ and $W_n\uparrow W$. 
Take $n$ with $\nu_q^\ast(W_n)>0$, where $\nu_q^\ast$
stands for the outer measure of $\nu_q$.

For any centered $\frac{1}{n}$-covering $\{\ball_j\}$ of $W_n$,
$$0<\nu_q^\ast(W_n)\leq \sum\nu_q^\ast(\ball_j)\leq \sum\nu_q(\ball_j)
\leq \sum\nu(\ball_j)^q r_j^{\underline{\tau}_\mu(q)-\varepsilon},$$
which implies
$$\overline{\mathscr{H}}_{\nu,\frac{1}{n}}^{q,\underline{\tau}_\mu(q)-
\varepsilon}(W_n)>0,$$
and
$$\mathscr{H}_\nu^{q,\underline{\tau}_\mu(q)-\varepsilon}(S_\nu)>0.$$

By definition this means $b_\nu(q)\geq\underline{\tau}_\mu(q)-\varepsilon$. 
Since $\varepsilon$ is arbitrary, we conclude that 
$b_\nu(q)\geq\underline{\tau}_\mu(q)$.
\qed

\subsection{$B_\nu(q)=\tau_\mu(q)$}
The following lemma gives a general estimate.

\begin{lemma}[see \cite{Bar}]\label{l4} 
Let $\mu$ be a probability measure on
$\partial\mathscr{A}^\ast$ and let $\nu=\gamma_\ast(\mu)$ be 
its projection onto $[0,1]$. One has $\tau_\nu(q)\leq\tau_\mu(q)$.
\end{lemma}

So it is sufficient to show that $B_\nu(q)\geq \tau_\mu(q)$. Following 
the spirit of Lemma \ref{l3}, we introduce
\begin{lemma}\label{l5}
Let $q\in \mathbb{R}$. For any $\varepsilon>0$, for $\nu_q$-almost 
every $x\in [0,1]$, for any $\delta>0$, there exists $r<\delta$, such that 
$$\nu_q(\ball(x,r))\leq \nu(\ball(x,r))^q r^{\tau_\mu(q)-\varepsilon}.$$
Then we have $B_\nu(q)\geq\tau_\mu(q)$.
\end{lemma}

\proof
The proof of the first assertion follows the same lines as the proof 
of the first assertion of Lemma \ref{l3}. So it is sufficient to only prove 
the second. One takes any family of $\{E_i\}$ such that $\cup E_i=S_\nu=[0,1]$ 
and for each $i$ one computes 
$\overline{\mathscr{P}}_\nu^{q,\tau_\mu(q)-\varepsilon}(E_i)$.

Denote by $N$ the $\nu_q$-null set, i.e. for any $x\in[0,1]\setminus N$, 
for any $\delta>0$, there exists $r<\delta$, such that $\nu_q(\ball(x,r))\leq 
\nu(\ball(x,r))^q r^{\tau_\mu(q)-\varepsilon}$.

By Besicovitch property, we can extract from $\{\ball(x,r)\}_{x\in E_i\setminus N}$ 
$C_\ball$ countable families $\{\ball_{j,k}\}_{1\leq j\leq C_\ball, \,k\geq 1}$ such
that $\cup_{j,k}\ball_{j,k}\supseteq E_i\setminus N$ and for any $j$,
$\{\ball_{j,k}\}_{k\geq 1}$ is a $\delta$-packing of $E_i\setminus N$.

Then one gets
$$\nu_q^\ast(E_i\setminus N)\leq \sum_{j,k}\nu_q^\ast(\ball_{j,k})\leq
\sum_{j,k}\nu_q(\ball_{j,k})\leq \sum_{j,k}  \nu(\ball_{j,k})^q
r_{j,k}^{\tau_\mu(q)-\varepsilon}.$$

So there exists $j$ such that $\sum_k  \nu(\ball_{j,k})^q
r_{j,k}^{\tau_\mu(q)-\varepsilon} \geq
\frac{1}{C_\ball}\nu_q^\ast(E_i\setminus N)$. Thus
$$\overline{\mathscr{P}}_\nu^{q,\tau_\mu(q)-\varepsilon}(E_i\setminus N)\geq
\frac{1}{C_\ball}\nu_q^\ast(E_i\setminus N),$$
which implies
$$\overline{\mathscr{P}}_\nu^{q,\tau_\mu(q)-\varepsilon}(E_i)\geq
\frac{1}{C_\ball}\nu_q^\ast(E_i),$$
and
$$\sum_i\overline{\mathscr{P}}_\nu^{q,\tau_\mu(q)-\varepsilon}(E_i)\geq
\frac{1}{C_\ball}\sum_i \nu_q^\ast(E_i)\geq
\frac{1}{C_\ball}\nu_q^\ast(S_\nu),$$
yielding that
$$\mathscr{P}_\nu^{q,\tau_\mu(q)-\varepsilon}(S_\nu)\geq
\frac{1}{C_\ball}\nu_q^\ast(S_\nu)>0.$$

Again by definition this means $B_\nu(q)\geq \tau_\mu(q)-\varepsilon$. 
Since $\varepsilon$ is arbitrary, $B_\nu(q)\geq\tau_\mu(q)$.
\qed
\medskip

So together with Lemma \ref{l4}, we conclude that $B_\nu(q)=\tau_\mu(q)$. 
And the proof of Theorem \ref{t2} is complete.

\section{The proof of Theorem \ref{t3}}\label{sec5}
We already know in \cite{Ben2, She} that for a certain range of $\alpha$, 
the Hausdorff dimension of the set $X_\mu(\alpha)$ is given by the value 
of the Legendre transform of~$b_\mu$ at~$\alpha$ whereas its packing 
dimension is the value of the Legendre transform of~$B_\mu$ at~$\alpha$. 
This means that the measure $\mu$ satisfies an extended multifractal 
formalism at $\alpha$. In this section, we show that for the same range 
of $\alpha$, the image measure $\nu$ also has this property. So the use 
of a Gray code is not necessary although we are dealing with non-doubling 
measures. As previously, we still set $c_1=c_2=2$.

For $\alpha\in(-\log s_1,-\log s_2)$, one finds $q_a,q_b$ such that 
$$-\theta'_a(q_a)=-\theta'_b(q_b)=\alpha.$$
Then one defines a new probability measure $\tilde{\mu}_q$ on the
symbolic space just as $\mu$ by replacing $\{a_1,a_2,b_1,b_2\}$ with
$\{\tilde{a}_1,\tilde{a}_2,\tilde{b}_1,\tilde{b}_2\}$, where
\begin{equation}\label{new para}
\tilde{a}_1=\frac{a_1^{q_a}}{a_1^{q_a}+a_2^{q_a}},\; 
\tilde{b}_1=\frac{b_1^{q_b}}{b_1^{q_b}+b_2^{q_b}}.
\end{equation}
And $\tilde{a}_2=1-\tilde{a}_1$, $\tilde{b}_2=1-\tilde{b}_1$. We also denote 
$\tilde{\nu}_q=\gamma_\ast(\tilde{\mu}_q)$.

Fix $\lambda< 1$ and define
\begin{multline*}
\varphi(t)=\limsup_{\delta\rightarrow 0}\frac{-1}{\log
  \delta}\log\sup\left\{\sum_i
\mu(\ball_i)^t\tilde{\mu}_q(\ball_i)\ :\ \right.\\ \left.(\ball_i)_i \textrm{ packing of
}S_\mu\ \textrm{with }\lambda\delta<r_i\leq\delta\vphantom
   {\sum_i}\right\},
\end{multline*}
then it is easy to compute
$$\varphi(t)=\log_2 \max\{a_1^t\tilde{a}_1+a_2^t\tilde{a}_2,\,b_1^t\tilde{b}_1+
b_2^t\tilde{b}_2\}.$$
So $\varphi(0)=0$. And the method of choosing $\{\tilde{a}_1,\tilde{b}_1\}$ in 
Formula~\eqref{new para} insures that $\varphi'(0)$ exists. In fact,
$$-\varphi'(0)=-\tilde{a}_1\log_2 a_1-\tilde{a}_2\log_2 a_2=
-\tilde{b}_1\log_2 b_1-\tilde{b}_2\log_2 b_2=\alpha.$$

Thus by \cite{Ben2,Pey} the set $X_\mu(\alpha)$ has full 
$\tilde{\mu}_q$-mass, and this implies that the set $X_\nu(\alpha)$ 
has full $\tilde{\nu}_q$-mass. Actually we have
\begin{lemma}\label{net}
For $\tilde{\nu}_q$-almost every $x\in [0,1]$, 
$$\liminf_{r\rightarrow 0} \frac{\log \nu(\ball(x,r))}{\log r}=
\liminf_{n\rightarrow \infty} \frac{\log \nu(I_n(x))}{\log 2^{-n}},$$
$$\limsup_{r\rightarrow 0} \frac{\log \nu(\ball(x,r))}{\log r}=
\limsup_{n\rightarrow \infty} \frac{\log \nu(I_n(x))}{\log 2^{-n}}.$$
\end{lemma}

\proof
First we revisit Lemma \ref{key l2}. We 
can obtain a constant $C_0>1$ and a subset $N\subseteq [0,1]$ 
such that $\tilde{\nu}_q(N)=0$, and that for any $x\in [0,1]\setminus N$, 
when $n$ is large enough,
\begin{equation}\label{similar}
C_0^{-(\sqrt{n}+1)}\leq \frac{\nu(I_n(x))}{\nu(I_n(x)^\pm)}\leq C_0^{\sqrt{n}+1}.
\end{equation}

For any $x\in [0,1]\setminus N$, when $r$ is small enough, we can 
choose $n$ large enough such that $2^{-(n+2)}<r\leq 2^{-(n+1)}$. 
So there exist at most two basic intervals $I_1, I_2\in\mathscr {F}_n$ 
such that $\ball(x,r)\subseteq I_1\cup I_2$; on the other hand, 
$\ball(x,r)$ must contain at least one basic interval of $(n+2)$-th generation, 
denoted by $I_3$. It is no restriction to assume that $x\in I_1$ and $x\in I_3$. 

Since $n$ is large, by \eqref{similar} we have
$$\nu(\ball(x,r))\leq \nu(I_1)+\nu(I_2)\leq (1+C_0^{\sqrt{n}+1})\nu(I_1).$$
And since $I_3\subseteq I_1$, 
$$\nu(\ball(x,r))\geq \nu(I_3)\geq C'_0\nu(I_1),$$
where $C'_0=\min\{a_i a_j,\,a_i b_j,\,b_i b_j :\,i,j=1,2 \}$.

To conclude, one notices that $I_1=I_n(x)$, and obtains 
$$\liminf_{r\rightarrow 0} \frac{\log \nu(\ball(x,r))}{\log r}=
\liminf_{n\rightarrow \infty} \frac{\log \nu(I_n(x))}{\log 2^{-n}}.$$
And the same goes for the upper limit.
\qed
\medskip

So we obtain from this lemma that $X_\nu(\alpha)\setminus N=
\gamma(X_\mu(\alpha))\setminus N$, while $\tilde{\nu}_q(N)=0$. 
Thus $\tilde{\nu}_q(X_\nu(\alpha))=1$, 
which yields
\begin{lemma}[see \cite{Ben2}]\label{minoration}  
Let $E=X_\nu(\alpha)$, then one has
$$\dim E\geq \esssup_{x\in E,\,\tilde{\nu}_q} \liminf_{r\rightarrow
0} \frac{\log \tilde{\nu}_q(\ball(x,r))}{\log r},$$
$$\Dim E\geq \esssup_{x\in E,\,\tilde{\nu}_q} \limsup_{r\rightarrow
0} \frac{\log \tilde{\nu}_q(\ball(x,r))}{\log r}.$$
\end{lemma}
\bigskip

However, using the same method as Lemma \ref{net}, the local H\"older 
exponent of~$\tilde{\nu}_q$ can be computed by applying general balls as 
well as dyadic intervals. For $\tilde{\nu}_q$-almost every $x\in [0,1]$, 
$$\liminf_{r\rightarrow 0} \frac{\log \tilde{\nu}_q(\ball(x,r))}{\log
  r}=\liminf_{n\rightarrow \infty} \frac{\log
  \tilde{\nu}_q(I_n(x))}{\log 2^{-n}},$$
$$\limsup_{r\rightarrow 0} \frac{\log \tilde{\nu}_q(\ball(x,r))}{\log
  r}=\limsup_{n\rightarrow \infty} \frac{\log
  \tilde{\nu}_q(I_n(x))}{\log 2^{-n}}.$$

Since dyadic intervals correspond to cylinders, we refer to \cite{Ben2,She} 
and present the proof of Theorem \ref{t3}.

\proof[Proof of Theorem \ref{t3}]
The strong law of large numbers shows that
$$\liminf_{n\rightarrow \infty} \frac{\log_2
  \tilde{\nu}_q(I_n(x))}{-n}=\min\{h(\tilde{a}),h(\tilde{b})\},$$
$$\limsup_{n\rightarrow \infty} \frac{\log_2
  \tilde{\nu}_q(I_n(x))}{-n}=\max\{h(\tilde{a}),h(\tilde{b})\},$$
for $\tilde{\nu}_q$-almost every $x$, where
$$h(\tilde{a})=-\sum_{i=1}^2 \tilde{a}_i\log_2 \tilde{a}_i\,\textrm{ and }\,
h(\tilde{b})=-\sum_{i=1}^2 \tilde{b}_i\log_2 \tilde{b}_i.$$

So it deduces from Lemma~\ref{minoration} that
$$\dim X_\nu(\alpha)\geq \min\{h(\tilde{a}),h(\tilde{b})\},$$
$$\Dim X_\nu(\alpha)\geq \max\{h(\tilde{a}),h(\tilde{b})\}.$$
And these two inequalities remain valid if we replace $\nu$ with $\mu$.

At the same time, one obtains
$$
\begin{array}{rccccll}
\medskip
h(\tilde{a}) &=& \theta_a(q_a)-q_a\theta'_a(q_a) &=& 
\theta^\ast_a(-\theta'_a(q_a)) &=& \theta^\ast_a(\alpha),\\
h(\tilde{b}) &=& \theta_b(q_b)-q_b\theta'_b(q_b) &=&
\theta^\ast_b(-\theta'_b(q_b)) &=& \theta^\ast_b(\alpha).
\end{array}
$$

Recall that the upper bounds of the dimensions of the level sets have
been given by Olsen \cite{Ols}. So for any $\alpha\in(-\log s_1,-\log
s_2)$, we have
$$\dim X_\mu(\alpha)=\dim
X_\nu(\alpha)=b_\mu^\ast(\alpha)=b_\nu^\ast(\alpha).$$
And for $\alpha\in(-\log s_1,-\log s_2)$ such that 
$$\max\{\theta^\ast_a(\alpha),\,\theta^\ast_b(\alpha)\}=B^\ast_\mu(\alpha),$$
we have
$$\Dim X_\mu(\alpha)=\Dim X_\nu(\alpha)=B_\mu^\ast(\alpha)=B_\nu^\ast(\alpha).$$
\qed

\begin{cor}
For any $q$ such that $B'_\mu(q)$ exists, denote $\alpha=-B'_\mu(q)$. 
If $\alpha\in(-\log s_1,-\log s_2)$, then we have 
$$\dim X_\mu(\alpha)=\dim X_\nu(\alpha)=b_\mu^\ast(\alpha)=b_\nu^\ast(\alpha),$$
$$\Dim X_\mu(\alpha)=\Dim X_\nu(\alpha)=B_\mu^\ast(\alpha)=B_\nu^\ast(\alpha).$$
\end{cor}

\proof
It is easy to see that when $B'_\mu(q)$ exists, $B'_\mu(q)$ coincides with 
either $\theta'_a(q)$ or $\theta'_b(q)$. Without loss of generality, we 
may assume that for such~$q$, $B'_\mu(q)=\theta'_a(q)$, then 
$B_\mu(q)=\theta_a(q)$. So 
$$\theta^\ast_a(-\theta'_a(q))=\theta_a(q)-q\theta'_a(q)=B_\mu(q)-qB'_\mu(q)=
B^\ast_\mu(-B'_\mu(q)),$$
which implies
$$\theta^\ast_a(\alpha)=B^\ast_\mu(\alpha).$$
But this means
$$\max\{\theta^\ast_a(\alpha),\,\theta^\ast_b(\alpha)\}=B^\ast_\mu(\alpha).$$
\qed

\begin{remark}\label{re}
One can before projection compose with an isometry of the symbolic 
space in Lemma \ref{key l2}, and thus in Theorem \ref{t2} and 
Theorem~\ref{t3}. To be precise, let $g:(\partial\mathscr{A}^\ast,d)
\rightarrow (\partial\mathscr{A}^\ast,d)$ be an isometry, and denote 
by $\nu_g$ the image measure of $\mu$ under $\gamma\circ g$. 
Then it is easy to see that for any two words $x$ and $y$, 
$g(x\wedge y)=g(x)\wedge g(y)$. So the proof of Inequality 
\eqref{key} is valid if we replace the measure $\nu$ with $\nu_g$.

Of course, Gray codes are isometries. As seen in \cite{Ben,Ben2,She}, 
if $g$ is a Gray code, then the measure $\nu_g$ becomes a doubling 
measure on $[0,1]$. But for general $g$, $\nu_g$ needs not be doubling.
\end{remark}


\begin{flushleft}
\bf Acknowledgements
\end{flushleft}
The author is grateful to Professor Jacques Peyri{\`e}re for his
patient guidance and helpful comments. The author would also like to
thank Mr.~Zhihui Yuan for the beneficial discussions the author had
with him.



\end{document}